\documentclass[11pt]{amsart}
\usepackage{amssymb,latexsym}
\usepackage{graphicx}
\newtheorem{thm}{Theorem}[section]
\newtheorem{lem}[thm]{Lemma}
\makeatletter\@addtoreset{equation}{section}\makeatother
\newcommand{\E}{{\mathbf E}}
\newcommand{\D}{{\mathbf D}}
\newcommand{\C}[1]{1-C^{#1}}
\newcommand{\CC}[1]{(1-C^{#1})}
\newcommand{\B}{{\mathbf B}}
\newcommand{\G}[2]{G_{#1}^{(#2)}}
\renewcommand{\H}[2]{H_{#1}^{(#2)}}
\newcommand\sumz[1]{\sum_{#1=0}^\infty}
\newcommand\superscript{$^1$}

\begin{document}

\title
[A Combinatorial Interpretation of The Numbers $6\,(2n)!\,/n!\,(n+2)!$]
{A Combinatorial Interpretation\\
of The Numbers $6\,(2n)!\,/n!\,(n+2)!$}
\author{Ira M. Gessel\superscript}
\address{Department of Mathematics\\
Brandeis University\\
Waltham MA 02454-9110}
\email{gessel@brandeis.edu}
\author{Guoce Xin}
\address{Department of Mathematics\\
Brandeis University\\
Waltham MA 02454-9110}
\email{maxima@brandeis.edu}
\date{August 5, 2004}
\thanks{$^1$Partially supported by NSF grant DMS-0200596}
\begin{abstract} It is well known that the numbers
$(2m)!\,(2n)!/m!\,n!\,(m+n)!$ are integers, but in general there is no known combinatorial interpretation
for them.
 When $m=0$ these numbers are the middle binomial coefficients
$\binom{2n}{n}$, and  when
$m=1$ they are twice  the Catalan numbers.  In this paper, we  give  combinatorial
interpretations for these numbers when $m=2$ or $3$.
\end{abstract}
\maketitle

\let\superscript=\relax

\section{Introduction}
The Catalan numbers
$$C_n=\frac{1}{n+1}\binom{2n}{n} =\frac{(2n)!}{n!\,(n+1)!}$$
are well-known integers that arise in many combinatorial problems.
Stanley \cite[pp. 219--229]{EC2}, gives
$66$ combinatorial interpretations of these numbers.

In 1874 E. Catalan \cite{catalan} observed that the numbers
\begin{equation}\frac{(2m)!\,(2n)!}{m!\,n!\,(m+n)!}\label{e-1}
\end{equation}
are integers, and their
number-theoretic
properties were  studied by several authors (see Dickson \cite[pp.
265--266]{dickson}).  For $m=0$,
\eqref{e-1} is the middle binomial coefficient
$\binom{2n}n$, and for $m=1$ it is $2C_n$.

Except for $m=n=0$, these integers are even, and it is convenient for our
purposes to
divide them by 2, so we consider the numbers
\begin{equation}
T(m,n)=\frac{1}{2} \; \frac{(2m)!\,(2n)!}{m!\,n!\,(m+n)!}.
\end{equation}

Some properties of these numbers are given in \cite{gessel}, where they are called ``super Catalan
numbers".  An intriguing problem is to find a combinatorial interpretation to the super Catalan
numbers. The following identity
\cite[Equation (32)]{gessel},  together with the symmetry property $T(m,n)=T(n,m)$ and the initial value $T(0,0)=1$,
shows that
$T(m,n)$ is a positive integer for all $m$ and $n$.
\begin{align}
\sum_{n} 2^{p-2n} {p\choose 2n} T(m,n) =T(m,m+p), \quad p\geq 0. \label{e-8}
\end{align}
Formula \eqref{e-8} allows us to construct recursively a set of cardinality $T(m,n)$ but we have not found any
natural description of it for $m\ge 2$.
Shapiro \cite{shapiro} gave a combinatorial interpretation to \eqref{e-8} in the case $m=1$,
which is the Catalan number identity
\[\sum_{n} 2^{p-2n}\binom p{2n} C_n=C_{p+1}.\]
A similar interpretation works for the case $m=0$ of \eqref{e-8} (when multiplied by 2), which is the
identity
\[\sum_{n} 2^{p-2n}\binom p{2n}\binom{2n}{n}=\binom{2p}{p}.\]

Another intriguing formula for $T(m,n)$, which does not appear in \cite{gessel}, is
\begin{equation}
\label{e-mo}
1+\sum_{m,n=1}^\infty C_m C_n x^m y^n=\biggl(1-\sum_{m,n=1}^\infty T(m,n)x^my^n\biggr)^{-1}.
\end{equation}
Although \eqref{e-mo} suggests a combinatorial interpretation for $T(m,n)$ based on a decomposition of pairs of
objects counted by Catalan numbers, we have not found such an interpretation.

In this paper, we give a combinatorial interpretation for
$T(2,n)=6\,(2n)!\,/n!\,(n+2)!$ for $n\ge1$ and for $T(3,n)= 60\,(2n)!\,/n!\,(n+3)!$ for $n\ge2$.
The first few values of $T(m,n)$ for $m=2$ and $m=3$ are as follows:
\[
\offinterlineskip\vbox
{\halign{\ \hfil\strut$#$\hfil\ \vrule&&\hfil\ $#$\ \cr
m\backslash n&0&1&2&3&4&5&6&7&8&9&10\cr
\noalign{\hrule}
2&3&2&3&6&14&36&99&286&858&2652&8398\cr
3&10&5&6&10&20&45&110&286&780&2210&6460\cr
}}
\]
We show that $T(2,n)$ counts pairs of Dyck paths of total length $2n$ with heights differing by at most 1.
We give two proofs of this result, one combinatorial and one using generating functions.
The combinatorial proof is based on
the  easily checked formula
\begin{equation}\label{e-2}T(2,n)=4C_n-C_{n+1}
\end{equation}
which we interpret by inclusion-exclusion.
%We  construct an injection
%from a set of cardinality $C_{n+1}$ to another set of cardinality
%$4C_n$, and thereby get a combinatorial interpretation for $4C_n-C_{n+1}$.

Our interpretation for
$T(3,n)$ is more complicated, and involves pairs of Dyck paths with height restrictions. Although we have the
formula $T(3,n) = 16C_n -8C_{n+1}+C_{n+2}$ analogous to $\eqref{e-2}$, we have not found a combinatorial
interpretation to this formula, and our proof uses generating functions.

Interpretations of the number $T(2,n)$ in terms of trees, related to each other, but not, apparently, to our
interpretation,  have been found by  Schaeffer \cite{gilles}, and by  Pippenger and  Schleich
\cite[pp.~34]{nicholas}.

\section{The Main theorem}

All paths in this paper have steps $(1,1)$  and $(1,-1)$, which we call \emph{up steps} and \emph{down steps}.
A step from a point
$u$ to a point $v$ is denoted by $u\to v$.
The \emph{level}
of a point in a path is defined to be its
$y$-coordinate.
A \emph{Dyck path} of \emph{semilength} $n$
(or of length $2n$)
is a path that starts at $(0,0)$,  ends at $(2n,0)$,
and never goes below level $0$.
It is well-known that the number of Dyck paths of semilength $n$ equals the Catalan
number $C_n$.
The
\emph{height} of  a path $P$,
denoted by $h(P)$, is the highest level it reaches.

Every nonempty Dyck path $R$ can be factored uniquely as $UPDQ$, where $U$ is an up step, $D$ is a down step,
and $P$ and $Q$ are Dyck paths. Thus the map $P\mapsto (P,Q)$ is a bijection from nonempty Dyck paths to
pairs of Dyck paths. Let $\B_n$ be the set of pairs of Dyck paths $(P,Q)$ of total semilength
$n$. This bijection gives $ | \B_n | = C_{n+1}$, so by \eqref{e-2}, we have $T(2,n)= 4C_n - | \B_n | $.

%
%Let $c(x)$ be the generating function for the Catalan numbers, so that
%\[c(x)=\sum_{n=0}^{\infty} \frac{1}{n+1}{2n\choose n} x^n
%=\frac{1-\sqrt{1-4x}}{2x}.\]
%Then $c(x)$ satisfies the functional equation $c(x)=1+xc(x)^2$. From
%equation \eqref{e-2}, the generating function for
%$T(2,n)$ is
%\begin{equation}
%\label{e-T2}
%\sum_{n=0}^\infty T(2,n) x^n
%  =\sum_{n=0}^\infty (4C_n-C_{n+1})x^n
%  = 4c(x)-\frac{c(x)-1}{x} = 4c(x)-c(x)^2.
%\end{equation}
%Now if we let $\B_n$ be the set of pairs of Dyck paths $(P,Q)$ of total semilength
%$n$, then $ |  \B_n |  $ has generating function $c(x)^2$, and we have $T(2,n)= 4C_n
%- |  \B_n |  $.

Our interpretatoin for $T(2,n)$ is a consequence of  the following Lemma \ref{l-main}.
We give two proofs of this lemma, one combinatorial and the other
algebraic. The algebraic proof will be given in the next section.

\begin{lem}\label{l-main}
For $n\ge 1$, $C_n$ equals the number of pairs of Dyck paths
$(P,Q)$ of total
semilength $n$, with $P$ nonempty and $h(P)\leq h(Q)+1$.
\end{lem}
\begin{proof}
Let $\D_n$ be the set of Dyck paths of semilength $n$,  and let $\E_n$ be
the set of pairs of Dyck paths $(P,Q)$ of total semilength $n$, with $P$
nonempty and $h(P)\leq h(Q)+1$.

For a given pair $(P,Q)$ in $\E_n$, since $P$ is nonempty, the last step of
$P$ must be a down step, say, $u\to v$. By replacing $u\to v$ in $P$ with an
up step $u\to v'$, we get a path $F_1$. Now raising $Q$ by two levels, we get
a path $F_2$. Thus $F:=F_1F_2$ is a path that ends at level $2$
and never goes below level $0$.  The point $v'$ belongs to both
$F_1$ and $F_2$, but we treat it as a point only in $F_2$, even if $F_2$
is the empty path.
The condition that
$h(P)\le h(Q)+1$ yields
$h(F_1)<h(F_2)$, which implies that the highest point of $F$ must belong to
$F_2$. See Figure \ref{f-pic1} below.

\begin{figure}[h]
\begin{center}
\includegraphics[scale=.8]{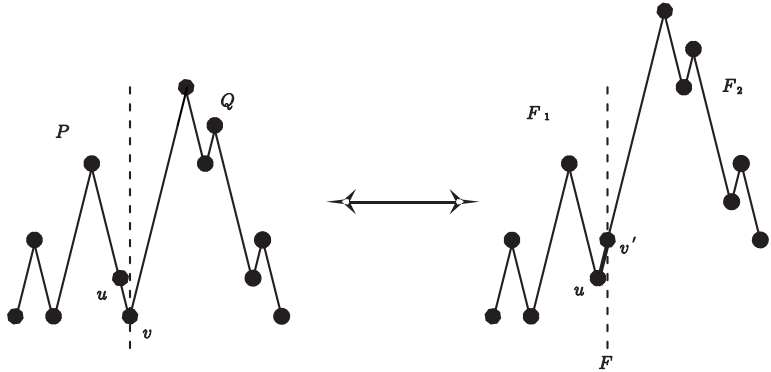}
\caption{First step of the bijection}
\label{f-pic1}
\end{center}
\end{figure}

Now
let $y$ be the leftmost highest point of $F$ (which is in $F_2$), and let
 $x\to y$ be the step in $F$ leading to $y$. Then $x\to y$ is  an up step.
By replacing $x\to y$ with a down step $x\to y'$, and lowering the part of
$F_2$ after $y$ by two levels,
we get a Dyck path $D\in \D_n$.  See Figure \ref{f-pic2} below.

\begin{figure}[h]
\begin{center}
\includegraphics[scale=.8]{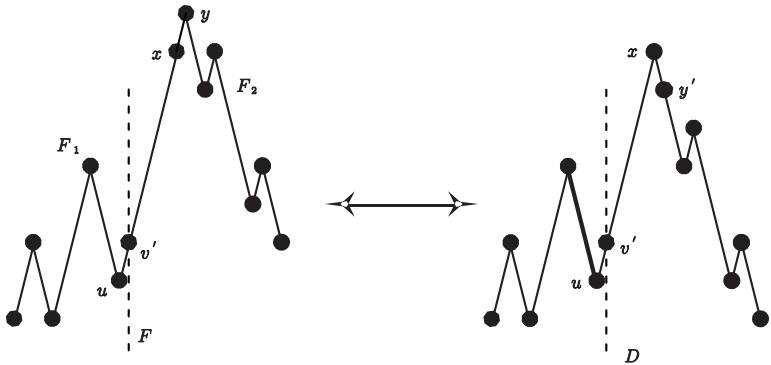}
\caption{Second step of the bijection}
\label{f-pic2}
\end{center}
\end{figure}

With the following two key observations, it is easy to see that the above procedure
gives a bijection from $\E_n$ to $\D_n$.
First,  $x$ in the final Dyck path $D$ is the rightmost
highest point. Second, $u$ in the intermediate path $F$ is the rightmost
point of level
$1$ in both $F$ and $F_1$.
\end{proof}

\begin{thm}\label{t-main}
For $n\ge 1$, the number $T(2,n)$ counts  pairs of Dyck paths $(P,Q)$ of total
semilength
$n$ with
$ |  h(P)-h(Q) |  \leq 1$.
\end{thm}
\begin{proof}
Let $F$ be the set of pairs of Dyck paths $(P,Q)$ with $h(P)\le h(Q)+1$, and
let $G$ be the set of pairs of Dyck paths $(P,Q)$ with $h(Q)\le h(P)+1$.
By symmetry, we see that
$ |  F |  = |  G |  $. Now we claim that the cardinality of $F$ is $2C_n$. This claim follows from
Lemma \ref{l-main} and the fact that
if $P$ is the empty path, then $h(P)\leq h(Q)+1$ for every $Q\in \D_n$.

Clearly we have that $F\cup G= \B_n$, and that $F\cap G$ is the set of
pairs of Dyck paths $(P,Q)$, with $ |  h(P)-h(Q) |  \leq 1$.
The theorem then follows from the following computation:
\begin{align*}
 |  F\cap G |    =  |  F |  + |  G |  - |  F\cup G |  =4C_n- |  \B_n |  =4C_n-C_{n+1}.
\end{align*}
%\vspace{-.7cm}
\end{proof}

\section{An Algebraic Proof and Further Results}

In this section we give an algebraic proof of Lemma \ref{l-main}.

Let $c(x)$ be the generating function for the Catalan numbers, so that
\[c(x)=\sum_{n=0}^{\infty} \frac{1}{n+1}{2n\choose n} x^n
=\frac{1-\sqrt{1-4x}}{2x}.\]
Then $c(x)$ satisfies the functional equation $c(x)=1+xc(x)^2$.
Let $C=xc(x)^2=c(x)-1$ and let $G_k$ be the generating function for Dyck paths
of height at most~$k$. Although $G_k$ is a rational function, a formula for $G_k$ in terms of $C$ will be of more
use to us than the explicit formula for $G_k$.
\begin{lem}
\label{l-G}
For  $k\ge -1$,
\begin{equation}
\label{e-G}
G_k=(1+C){1-C^{k+1}\over 1-C^{k+2}}.
\end{equation}
\end{lem}
\begin{proof}
Let $P$ be a path of height at most $k\ge1$. If $P$ is nonempty then
$P$ can be factored as $UP_1DP_2$, where $U$ is an up step, $P_1$ is a Dyck path of
height at most $k-1$ (shifted up one unit), $D$ is a down step, and $P_2$ is a Dyck path of height at most $k$. Thus
$G_k=1+ xG_{k-1}G_k$, so
$G_{k}=1/(1-xG_{k-1})$. Equation \eqref{e-G} clearly holds for $k=-1$ and $k=0$. Now
suppose that for some $k\ge1$,
\[G_{k-1}=(1+C){1-C^{k}\over 1-C^{k+1}}.\]
Then the recurrence, together with the formula $x=C/(1+C)^2$, gives
\begin{align*}G_k &=\left[1-x(1+C){1-C^{k}\over 1-C^{k+1}}\right]^{-1}
=\left[1-\frac C{1+C}\,\frac{1-C^{k}}{ 1-C^{k+1}}\right]^{-1}\\
&=\left[\frac{1-C^{k+2}}{(1+C)(1-C^{k+1})}\right]^{-1}
=(1+C)\frac{1-C^{k+1}}{1-C^{k+2}}.
\end{align*}

%\vspace{-.5cm}
\end{proof}

We can prove Lemma \ref{l-main}  by showing that
$\sum_{n=0}^\infty G_{n+1}(G_n-G_{n-1})=1+2C$; this is equivalent to the
statement that the
number of pairs $(P,Q)$ of Dyck paths of semilength $n>0$ with $h(P)\leq h(Q)+1$ is
$2C_n$.

\begin{thm}
\label{t-firstsum}
\[\sum_{n=0}^\infty(G_n-G_{n-1})G_{n+1}=1+2C.\]
\end{thm}

\begin{proof}
Let \[\Psi_k=\sum_{n=k}^\infty \frac{C^n}{1-C^{n}}.\]
Thus if $j<k$ then
\begin{equation}
\label{e-Psi}
\Psi_j=\sum_{n=j}^{k-1}\frac{C^n}{1-C^{n}}+\Psi_k.
\end{equation}

We have \[G_nG_{n+1}=(1+C)^2\frac{\C{n+1}}{\C{n+3}}\]
and
\[G_{n-1}G_{n+1}=(1+C)^2\frac{\CC n\CC{n+2}}{\CC{n+1}\CC{n+3}}.\]
Let
\[S_1=\sum_{n=0}^\infty \left(\frac{\C{n+1}}{\C{n+3}} -1\right)\]
and
\[S_2=\sum_{n=0}^\infty\left(\frac{\CC n\CC{n+2}}{\CC{n+1}\CC{n+3}}-1\right).\]
Then $\sum_{n=0}^\infty(G_n-G_{n-1})G_{n+1}=(1+C)^2(S_1-S_2)$.
We have
\[\frac{\C{n+1}}{\C{n+3}}-1=-\frac{\CC2C^{n+1}}{\C{n+3}},\]
so $S_1=-\CC2 C^{-2}\Psi_3$, and
\[\frac{\CC n\CC{n+2}}{\CC{n+1}\CC{n+3}}-1
=-\frac{(1-C)C^n}{(1+C)\CC{n+1}} -\frac{\CC3 C^{n+1}}{(1+C)\CC{n+3}},\]
so
\[S_2=-\frac{1-C}{1+C}C^{-1}\Psi_1 -\frac{1-C^3}{1+C}C^{-2}\Psi_3.\]
Therefore
\begin{align*}
S_1-S_2&=\frac{1-C}{1+C}C^{-1}\Psi_1
  +\left(\frac{1-C^3}{1+C}-(1-C^2)\right)C^{-2}\Psi_3\\
&=\frac{1-C}{1+C}C^{-1}(\Psi_1-\Psi_3)\\
&=\frac{1-C}{1+C}C^{-1}\left(\frac C{1-C}+\frac {C^2}{1-C^2}\right)
  =\frac{1+2C}{(1+C)^2}.
\end{align*}
Thus
$(1+C)^2(S_1-S_2)=1+2C$.
\end{proof}

By similar reasoning, we could prove Theorem \ref{t-main} directly: The generating function for pairs of paths with
heights differing by at most 1 is
\[\sum_{n=0}^\infty (G_n-G_{n-1})(G_{n+1}-G_{n-2}),\]
where we take $G_{-1}=G_{-2}=0$, and a calculation like that in the proof of Theorem \ref{t-firstsum} shows that
this is equal to

\begin{align*}
1+2C-C^2&=4c(x)-c(x)^2-2=4c(x) -\frac{c(x) -1}{x} -2\\
&=1+\sum_{n=1}^\infty (4C_n-C_{n+1})x^n
=1+\sum_{n=1}^\infty T(2,n)x^n.
\end{align*}

Although the fact that the series in Theorem \ref{t-firstsum}
telescopes may seem surprising, we shall see in Theorem \ref{t-general} that it is a special case of a
very general result on sums of generating functions for
Dyck paths with restricted heights.
%%%\begin{lem}
%%%\label{l-finite}
%%%If $\sum_{j=1}^s u_j=0$ then
%%%\[\sum_{j=1}^s u_j\Psi_{j}
%%%=\sum_{n=1}^{s-1}\biggl(\sum_{j=1}^{n}u_j\biggr)
%%%  \frac{C^n}{1-C^n}.\]
%%%\end{lem}
%%%
%%%\begin{proof}
%%%We have
%%%\[\sum_{j=1}^s u_j\Psi_{j}
%%%  =\sum_{j=1}^s u_j\sum_{n=j}^\infty\frac{C^n}{1-C^n}
%%%  =\sum_{n=1}^\infty \frac{C^n}{1-C^n}\sum_{j=1}^{\min(s,n)}u_j
%%%  =\sum_{n=1}^{s-1} \frac{C^n}{1-C^n}\sum_{j=1}^{n}u_j.
%%%\]
%%%\end{proof}
%%%\comment{Figure out how to move up the \kern -5pt\qed\hskip 0pt plus 1 filll.}

\begin{lem}\label{l-general}
Let $R(z,C)$ be a rational function of $z$ and $C$ of the form
\[\frac{zN(z,C)}{\prod_{i=1}^m(1-zC^{a_i})},\]
where $N(z,C)$ is a polynomial in $z$ of degree less than $m$, with coefficients
that are rational functions of $C$, and the $a_i$ are distinct positive integers.
Let $L=-\lim_{z\to\infty}R(z,C)$. Then
\[\sumz n R(C^n,C)=Q(C)+L\Psi_1,\]
where $Q(C)$ is a rational function of $C$.
\end{lem}

\begin{proof}
First we show that the lemma holds for $R(z,C)=z/(1-zC^a)$.
In this case,
$L=-\lim_{z\to\infty} R(z,C)=C^{-a}$ and
\[\sumz n R(C^n,C)=
\sumz n \frac{C^n}{1-C^{n+a}}=C^{-a}\sum_{n=a}^\infty \frac{C^{n}}{1-C^{n}}
=-\sum_{n=0}^{a-1}\frac{C^{n-a}}{1-C^n}+C^{-a}\Psi_1.\]

Now we consider the general case. Since $R(z,C)/z$ is a proper rational function of
$z$, it has a partial fraction expansion
\[\frac1z R(z,C)=\sum_{i=1}^m\frac{U_i(C)}{1-zC^{a_i}}\]
for some rational functions $U_i(C)$, so
\[R(z,C)=\sum_{i=1}^m U_i(C)\frac{z}{1-zC^{a_i}}.\]
The general theorem then follows by applying the special case to each summand.
\end{proof}

\begin{thm}\label{t-general}
If $i_1, i_1',\ldots, i_m$ are distinct integers, then
\[\sumz n(G_{n+i_1}-G_{n+i_1'})G_{n+i_2}G_{n+i_3}\cdots G_{n+i_m}\]
is a rational function of $C$.
\end{thm}
\begin{proof}
Apply Lemma \ref{l-general} to \[\sumz n\left(G_{n+i_1}G_{n+i_2}G_{n+i_3}\cdots
G_{n+i_m}-1\right)\]
and
\[\sumz n \left(G_{n+i_1'}G_{n+i_2}G_{n+i_3}\cdots G_{n+i_m}-1\right).\]

%\vspace{-5mm}
\end{proof}

\section{A combinatorial interpretation for $T(3,n)$}

It is natural to ask whether there are combinatorial interpretations to $T(m,n)$ for $m>2$ similar to
Theorem \ref{t-main}. It is straightforward (with the help of a computer algebra system) to evaluate sums like
\eqref{t-general} that count
$m$-tuples of paths with height restrictions, and we find that sums like that in Theorem \ref{t-general}
involving products of $m$ path generating functions may generally be expressed in the form
\[R_1(x)+R_2(x)\sum_{n=0}^\infty T(m,n) x^n,\]
where $R_1(x)$ and $R_2(x)$ are rational functions of $x$. However, for general $m$ we have not found a set of
$m$-tuples of paths counted by $T(m,n)$. We have found a set of paths counted by $T(3,n)$, though it is
not as simple as one would like.

We need to consider paths that end at levels greater than 0.
Let us define a \emph{ballot path} to be a path that starts at level 0
and never goes below level 0.
%, and define a \emph{generalized ballot path}
%to be a path that never goes below level 0.

In the previous section all our
paths had an even number of steps, so it was natural to assign a path with $n$ steps the weight $x^{n/2}$. We shall
continue to weight paths in this way, even though some of our paths now have odd lengths.

Let $G_k^{(j)}$ be the generating function for ballot paths of height at most $k$ that end at level $j$.

\begin{lem}
\label{t-Gij}
For $0\le j\le k+1$ we have
\begin{equation}
\label{e-Gj}
G_k^{(j)}=C^{j/2}(1+C)\frac{\C{k-j+1}}{1-C^{k+2}}
\end{equation}

\end{lem}

\begin{proof}
The case $j=0$ is Lemma \ref{l-G}. Now let $W$ be a ballot
path counted by $G_k^{(j)}$, where $j>0$, so that $W$ is of height at most $k$ and $W$
ends at level $j$.  Then
$W$ can be factored uniquely as $W_1UW_2$, where $W_1$ is a path of height at most $k$ that
ends at level 0 and $W_2$ is a path from level 1 to level $j$ that never goes above level
$k$ nor below level 1.
Using $\sqrt x=\sqrt{C/(1+C)^2}=\sqrt C/(1+C)$, we obtain
\[G_k^{(j)}=G_k\sqrt x G_{k-1}^{(j-1)}=(1+C)\frac{\C {k+1}}{\C{k+2}}
\cdot \frac{\sqrt C}{1+C}\cdot
G_{k-1}^{(j-1)}=\sqrt C\frac{\C {k+1}}{\C{k+2}}G_{k-1}^{(j-1)},\]
and \eqref{e-Gj} follows by induction.

\end{proof}

We note an alternative formula that avoids  half-integer powers of $C$,
\[G_k^{(j)}=x^{j/2}(1+C)^{j+1}\frac{\C{k-j+1}}{1-C^{k+2}},\]
which follow easily from \eqref{e-Gj}
% \eqref{e-Gij},
and the formula $\sqrt x=\sqrt C/(1+C)$.

Although we will not need it here, there is a similar formula for the generating function $G_k^{(i,j)}$
for paths of height at most $k$ that start at level $i$, end at level $j$, and never go below level 0:
for $0\le i\le j\le k+1$ we have
\begin{equation}
\label{e-Gij}
G_k^{(i,j)}=
C^{(j-i)/2}(1+C)\frac{\CC{i+1}(1-C^{k-j+1})}{\CC{}\CC{k+2}},
\end{equation}
with $G_k^{(i,j)}=G_k^{(j,i)}$ for $i>j$.

We note two variants of \eqref{e-Gij}, also valid for $0\le i\le j\le k+1$:
\begin{align*}
G_k^{(i,j)}&=x^{(j-i)/2}(1+C)^{j-i+1}\frac{\CC{i+1}(1-C^{k-j+1})}{\CC{}\CC{k+2}},\\
  &=x^{-1/2}C^{(j-i+1)/2}\frac{\CC{i+1}\CC{k-j+1}}{\CC{}\CC{k+2}}.
\end{align*}

It is well known that $G_k^{(i,j)}$ is $x^{(j-i)/2}$ times a rational function of $y$, and it is
useful to have an explicit formula for it as a quotient of polynomials. (See Sato and Cong \cite{sato-cong} and
Krattenthaler \cite{kratt}.) Let us define polynomials $p_n=p_n(x)$
 by
\[p_n(x)=\sum_{0\le k\le n/2}(-1)^k{n-k\choose k}x^k.\]
%It follows easily that
%\[\sum_{n=0}^\infty p_n z^n={1\over 1-z+xz^2},\]
%and a partial fraction expansion gives
%\begin{equation}
%\label{e-p}
%p_{n}
% ={1\over\sqrt{1-4x}}\left[\left(1+\sqrt{1-4x}\over 2\right)^{n+1}
%   -\left(1-\sqrt{1-4x}\over 2\right)^{n+1}\right].
%\end{equation}
The first few values are
\begin{align*}
p_0&=1\\
p_1&=1\\
p_2&=1-x\\
p_3&=1-2x\\
p_4&=1-3x+x^2\\
p_5&=1-4x+3x^2\\
p_6&=1-5x+6x^2-x^3\\
\end{align*}

These polynomials can be expressed in terms of the
Chebyshev polynomials of the second kind $U_n(x)$  by
%\[\sum_{n=0}^\infty U_n(x) z^n={1\over 1-2xt+t^2},\]
%so
\[p_n(x)=x^{n/2}U_n\left(1\over 2\sqrt x\right).\]
%It follows from \eqref{e-p} that
It is not difficult to show that
\[p_n=\frac{\C{n+1}}{(1-C)(1+C)^n},\]
and thus we obtain
\[G_k^{(i,j)}=x^{(j-i)/2}\frac{p_i p_{k-j} }{p_{k+1}},\]
for $0\le i\le j\le k$, and in particular,
$G_k^{(j)}=x^{j/2} p_{k-j}/p_{k+1}$ and $G_k=p_k/p_{k+1}.$

We can now describe our combinatorial interpretation of $T(3,n)$: $T(3,n)$ counts triples of ballot paths whose
heights are  $k$, $k-2$, and $k-4$ for some $k$, ending at levels 4, 3, and 2, together with some additional paths
of height at most 5. (Note that if a path of height $k-4$ ends at level 2, then $k$ must be at least 6.) More
precisely, let
$\H kj$ be the generating function for ballot paths of height
$k$ that end at level
$j$, so that $\H kj = \G kj - \G{k-1}j$.
Then we have:

\begin{thm}
\begin{equation}
\label{e-T3}
1+\sum_{n=0}^\infty T(3,n+1) x^n = \sqrt x\sum_{k=6}^\infty \H k4 \H{k-2}3 \H{k-4}2
+2G_1+2G_2+G_3+G_5.
\end{equation}
\end{thm}

\begin{proof} First note that $T(m,n)= \tfrac 12 (-1)^n 4^{m+n}\binom{m-\tfrac12}{m+n}$, so
\begin{align*}
\tfrac12 (-1)^m (1-4x)^{m-1/2}&=\sum_{n=-m}^\infty \tfrac 12 (-1)^n 4^{m+n}\binom{m-\tfrac12}{m+n}x^{m+n}\\
   &=\sum_{n=-m}^0 \tfrac 12 (-1)^n 4^{m+n}\binom{m-\tfrac12}{m+n}x^{m+n}
    +\sum_{n=1}^\infty T(m,n) x^{m+n}.
\end{align*}
Setting $m=3$, dividing both sides by $x^4$, and subtracting the first sum on the right from both sides gives
\begin{equation}
\label{e-5/2}
-\frac{(1-4x)^{5/2}}{2x^4} - \frac{10}x +\frac{15}{x^2}-\frac{5}{x^3}+\frac1{2x^4}
=\sum_{n=0}^\infty T(3,n+1)x^n.
\end{equation}

Using the method described in Lemma \ref{l-general} and Theorem \ref{t-general}, we find, with the
help of Maple, that the sum
\[ \sqrt x\sum_{k=6}^\infty \H k4 \H{k-2}3 \H{k-4}2\]
is equal to
\[-\frac{(1-4x)^{5/2}}{2x^4} - \frac{10}x +\frac{15}{x^2}-\frac{5}{x^3}+\frac1{2x^4} +1 -\frac2{1-x}
-2\frac{1-x}{1-2x} -\frac{1-2x}{1-3x+x^2}-\frac{1-4x+3x^2}{1-5x+6x^2-x^3}.
\]
Then \eqref{e-T3} follows from \eqref{e-5/2}, the formula $G_k=p_k/p_{k+1}$, and the formulas for $p_k$, $k=1,\dots,
6.$
\end{proof}

\end{document}